\documentclass{elsart}
\usepackage{amssymb,amsmath}

\newtheorem{theorem}{Theorem}[section]
\newtheorem{lemma}{Lemma}[section]
\newtheorem{proposition}{Proposition}[section]
\newtheorem{remark}{Remark}[section]

\newtheorem{definition}{Definition}[section]

\makeatletter\def\theequation{\arabic{section}. \arabic{equation}}\makeatother

\begin{document}

\begin{frontmatter}

\title{\bf Long-time extinction of solutions of some semilinear parabolic equations}

\author[label1,label2]{Yves Belaud, Andrey Shishkov}

\address[label1]{{\it Laboratoire de Math\'ematiques et Physique Th\'eorique (UMR
CNRS 6083),}\\{\it F\'ed\'eration Denis Poisson,} \it{Universit\'e
Fran\c{c}ois Rabelais,}\\{\it Parc de Grandmont,} \it{37200 Tours,
France}}

\address[label2]{{\it Institute of Applied Mathematics and Mechanics of NAS of Ukraine},\\
{\it R. Luxemburg str. 74, 83114 Donetsk, Ukraine}}

\ead{yves.belaud@univ-tours.fr,shishkov@iamm.ac.donetsk.ua}

\begin{abstract}
We study the long time behaviour of solutions of semi-linear parabolic equation of
the following type $\partial_t u-\Delta u+a_0(x)u^q=0$ where
$a_0(x) \geq d_0 \exp\left(-\frac{\omega(|x|)}{|x|^2}\right)$, $d_0>0$, $1>q>0$ and $\omega$
a positive continuous radial function. We give a Dini-like condition on the function $\omega$
by two different method which implies that any solution of the above equation vanishes in a finite time.
The first one is a variant of a local energy method and the second one is derived from semi-classical
limits of some Schr\"odinger operators.
\end{abstract}

\end{frontmatter}

\date{\today}

\setcounter{section}{0}

\section{Introduction}

\def\theequation{1.\arabic{equation}}\makeatother

\setcounter{equation}{0}

Let $\Omega\subset\mathbb{R}^N,\ N\geq1,$ be a bounded domain with
$C^1$-boundary, $0\in\Omega$. The aim of this paper is to
investigate the time vanishing properties of generalized (energy)
solutions of initial-boundary problem to a wide class of
quasilinear parabolic equations with the model representative:
\begin{eqnarray} \label{parabolicequation}
\left\{ \begin{array}{rllcl} u_t- \Delta u + a_0(x) |u|^{q-1}u & =
& 0 &
\rm in & \Omega\times(0,\infty) , \\
\frac{\partial u}{\partial n} & = & 0 & \rm on & \partial\Omega\times(0,\infty), \\
u(x,0) & = & u_0(x) & \rm on & \Omega,
\end{array} \right.
\end{eqnarray}
where $0<q<1,\ a_0(x)\geq0$ and $u_0 \in L_2(\Omega)$. It is easy
to see that if $a_0(x)\geq\varepsilon>0$, then the comparison with
the solution of corresponding ordinary equation
$\varphi_t+\varepsilon|\varphi|^{q-1}\varphi=0$ implies that the
solution $u(x,t)$ of \eqref{parabolicequation} vanishes for $t\geq
T_0=\varepsilon^{-1}(1-q)^{-1}\|u_0\|_{L_\infty}^{1-q}$. The
property that any solution of problem \eqref{parabolicequation}
becomes identically zero for $t$ large enough is called the time
compact support property (TCS-property). On the opposite, if
$a_0(x)\equiv0$ for any $x$ from some connected open subset
$\omega\subset\Omega$, then any solution $u(x,t)$ of problem
\eqref{parabolicequation} is bounded from below by
$\sigma\exp(-t\lambda_\omega)\varphi_\omega(x)$ on
$\omega\times(0,\infty)$, where $\sigma=\text{ess}\inf_\omega
u_0>0$, $\lambda_\omega$ and $\varphi_\omega$ are first eigenvalue
and corresponding eigenfunction of $-\Delta$ in
$W_0^{1,2}(\omega)$. It was Kondratiev and Veron \cite{KV97} who
first proposed  a method of investigation of conditions of
appearance of TCS-property in the case of general potential
$a_0\geq0$. They introduced the fundamental states of the
associated Schr\"odinger operator
\begin{equation}\label{1.1*}
\mu_n=\inf\bigg\{ \int_\Omega(|\nabla\psi|^2+2^na_0(x)\psi^2)\,dx\
:\ \psi\in W^{1,2}(\Omega),\quad \int_\Omega\psi^2dx=1
\bigg\},\quad n\in\mathbb{N},
\end{equation}
and proved that, if
\begin{equation}\label{1.1**}
\sum_{n=0}^\infty\mu_n^{-1}\ln (\mu_n)<\infty,
\end{equation} then
\eqref{parabolicequation} possesses the TCS-property. Starting
from  condition \eqref{1.1**} in \cite{BHV01} an explicit
conditions of appearance of TCS-property in terms of potential
$a_0(x)$ was obtained. The analysis in \cite{BHV01} was based on
the so-called semiclassical analysis \cite{He89}, which uses sharp
estimates of the spectrum of the Schr\"odinger operator
\cite{Cw77,LT76,Ro72}. Particularly, in the case of existence of
the radially symmetric minorant
\begin{equation}\label{1.4}
a_0(x)\geq
d_0\exp\Big(-\frac{\omega(|x|)}{|x|^2}\Big):=a(|x|)\quad\forall\,x\in
\Omega,\ d_0>0
\end{equation}
the following statements was
obtained in \cite{BHV01}:
\begin{proposition} [Th. 4.5 from \cite{BHV01}] In equation \eqref{parabolicequation}
let $a_0(x)=a(|x|)$, where $a(r)$ is defined by
 $\eqref{1.4}$. Let $u_0(x)\geq\nu>0\
 \forall\,x\subseteq\overline\Omega$ and $\omega(r)\to\infty$ as
 $r\to0$. Then arbitrary solution $u$ of problem
 \eqref{parabolicequation} never vanishes on $\Omega$.
\end{proposition}

\begin{proposition}[Corollary of Th. 3.1 in \cite{BHV01}] If in assumption
\eqref{1.4}\\ $a_0(x)= a(|x|)$ and $\omega(r)= r^\alpha$ with
$0 < \alpha < 2$ then an arbitrary solution of
\eqref{parabolicequation} enjoys the TCS-property.
\end{proposition}

Thus, an open problem is to find sharp border which distinguish
two different decay properties of solutions, described in
Proposition 1.1 and Proposition 1.2. Moreover, the method of
investigations used in \cite{BHV01,KV97} exploits essentially some
regularity properties of solutions under consideration,
particularly, sharp upper estimates of
$\|u(x,t)\|_{L_\infty(\Omega)}$ with respect to $t$. Such an
estimate is difficult to obtain or is unknown for solutions of
equations of more general structure than
\eqref{parabolicequation}. Particularly, it is absolutely
impossible to have any information about such a behaviour for
higher order parabolic equations. We propose here some new energy
method of investigations, which deals with energy norms of
solutions $u(x,t)$ only and, therefore, may be applied,
particularly, for higher order equations, too.

We suppose that function $\omega(s)$ from condition \eqref{1.4}
satisfies the conditions:
\begin{itemize}
\item [$(A_1)$] $\quad \omega(r)$ is continuous and nondecreasing
function $\forall \, r \geq 0$, \item [$(A_2)$] $\quad
\omega(0)=0,\ \omega(r)>0\ \forall \, r>0. $ \item [$(A_3)$]
$\quad \omega(s) \leq \omega_0< \infty \quad\forall \,s \in
\mathbb{R}^1_+$
\end{itemize}

Our main result reads as follows
\begin{theorem}\label{MR} Let $u_0(x)$ be an arbitrary function from $L_2(B_1)$, let
function $\omega(r)$ from  $\eqref{1.4}$ satisfy assumptions
$(A_1),\ (A_2),\ (A_3)$ and the following main condition:
\begin{equation}\label{1.2}
 \int\limits_0^c
\frac{\omega(s)}s\,ds<\infty \quad(\text{Dini like condition}).
\end{equation}
Suppose also that  $\omega(r)$ satisfies the following technical
condition
\begin{equation}\label{1.5}
\frac{s\omega'(s)}{\omega(s)} \leq 2-\delta\quad\forall\,s\in(0,s_0),\
s_0>0,\ 2>\delta>0.
\end{equation}
 Then an arbitrary energy solution $u(x,t)$ of the problem
$(\ref{parabolicequation})$ vanishes on $\Omega$ in some finite
time $T<\infty$.
\end{theorem}
In the sequel of the paper we show that the sufficiency of the
Dini condition \eqref{1.2} for the validity of TCS-property can be
proved also by the methods from \cite{KV97,BHV01} if one uses
$L_\infty$ estimates of solution $u(x,t)$ of problem
\eqref{parabolicequation}. This leads to the following result.

\begin{proposition}
The assertion of Theorem~1.1 holds if the function $\omega(s)$
satisfies conditions $(A_1)$--$(A_3)$, the Dini condition \eqref{1.2}
and the following similar to \eqref{1.5} technical conditions:
\begin{equation}\label{1.6}
\omega(s)\geq s^{2-\delta}\quad\forall\,s\in(0,s_0),\ s_0>0,\
2>\delta>0,
\end{equation}
\begin{equation}\label{1.7}
\text{the function $\frac{\omega(s)}{s^2}$ is decreasing on
}(0,s_0).
\end{equation}
\end{proposition}
\begin{remark}
It is easy to check that the function $\omega(s)=(\ln
s^{-1})^{-\beta}$ satisfies all the conditions of Theorem~1.1 and
Proposition~1.3 for arbitrary $\beta>1$.
\end{remark}

\section{\bf The proof of main result}
\def\theequation{2.\arabic{equation}}\makeatother
\setcounter{equation}{0}

The proof of Theorem~\ref{MR} is based on some variant of the local
energy method, which was developed, particularly in \cite{Sh1,Sh2}.
First, we introduce the following families of subdomains:
\[
\Omega(\tau)=\Omega\cap\{ |x|>\tau \},\quad
Q_s^{(T)}(\tau)=\Omega(\tau)\times(s,T),\quad T<\infty.
\]
\begin{definition}\label{d1}
An energy solution of problem $\eqref{parabolicequation}$ is the
function $u(x,t)\in L_2(0,T; W_2^1(\Omega))$:
$\frac{\partial u}{\partial t}\in L_2(0,T;(W_2^1(\Omega))^*), \
u(x,0)=u_0$, satisfying the following integral identity:
\begin{equation}\label{2.1**}
\int_0^T\langle u_t,\varphi \rangle
dt+\int_{\Omega\times(0,T)}(\nabla_xu,\nabla_x\varphi)\,dxdt+\int_{\Omega\times(0,T)}
a_0(x)|u|^{q-1}u\varphi\,dxdt=0
\end{equation}
for arbitrary $\varphi\in L_2(0,T; W_2^1(\Omega))\ \forall\,
T<\infty$.
\end{definition}
\begin{lemma}\label{l:2.1}
An arbitrary energy solution $u$ of the problem
\eqref{parabolicequation} satisfies the following global a priory
estimate
\begin{multline}\label{2.1*}
\int_\Omega|u(x,\hat{t})|^2\,dx+
\int_{Q_0^{(\hat{t})}(0)}(|\nabla_xu|^2+a(|x|)|u|^{q+1})\,dxdt\\*
\leq\int_\Omega|u_0|^2\,dx:=y_0,\quad\forall\,\hat{t}>0.
\end{multline}
\end{lemma}
%\begin{proof}
 Testing integral identity \eqref{2.1**} by
$\varphi(x,t)=u(x,t)\xi(x)$, where $\xi(x)$ is arbitrary
$C^1$-function, due to formula of integration by parts
\cite{Lio}, we derive the following equality:
\begin{multline}\label{2.2**}
2^{-1}\int_{\Omega}|u(x,\hat t)|^2\xi\,dx+
\int_{\Omega\times(s,\hat
t)}(|\nabla_xu|^2\xi+(\nabla_xu,\nabla_x\xi)u)\,dxdt\\+
\int_{\Omega\times(s,\hat t
)}a_0\xi|u|^{q+1}\,dxdt=2^{-1}\int_\Omega |u(x,s)|^2\xi\,dx, \quad
0\leq s< \hat t<\infty.
\end{multline}
Let $\eta(r)\in C^1(\mathbb{R}^1)$ be such that
$0\leq\eta(r)\leq1\ \forall\,r\in\mathbb{R}^1, \ \eta(r)=0$ if
$r\leq0,\ \eta(r)=1$ if $r>1$. Fix arbitrary numbers $\tau>0,\
\nu>0$ and  test
\eqref{2.2**} by
$$
\xi(x)=\xi_{\tau,\nu}(|x|):=\eta\Big( \frac{|x|-\tau}\nu\Big).
$$
Then passing to the limit $\nu\to 0$ we obtain
\begin{multline}\label{2.1}
2^{-1}\int_{\Omega(\tau)}|u(x,\hat{t})|^2\,dx+
\int_s^{\hat{t}}\int_{\Omega(\tau)}(|\nabla_xu|^2+a_0(x)|u|^{q+1})\,dxdt\\=
2^{-1}\int_{\Omega(\tau)}|u(x,s)|^2\,dx+\int_s^{\hat{t}}\int_{|x|=\tau}u\frac{\partial
u}{\partial n}\,d\sigma dt \quad \forall\,
\hat{t}:s<\hat{t}<\infty.
\end{multline}
From \eqref{2.1} with $\tau=0,\ s=0$ the necessary global
estimate \eqref{2.1*} follows.
%\end{proof}
Further we will denote by $c,c_i$ different positive constants
which depend on known parameters of the problem
\eqref{parabolicequation} only. Let us introduce the energy functions
related to a fixed energy solution $u$ of problem
\eqref{parabolicequation}:
\begin{equation}\label{2.7*}
H({t},\tau)=\int_{\Omega(\tau)}|u(x,{t})|^2\,dx,\qquad
I_s^{(v)}(\tau)=\int_{Q_s^{(v)}(\tau)}(|\nabla_xu|^2+a(|x|)|u|^{q+1})\,dxdt
\end{equation}
\[
E(t,\tau)=\int_{\Omega(\tau)}(|\nabla_xu(x,t)|^2+a(|x|)|u(x,t)|^{q+1})\,dx,
\]
\[
J_s^{(v)}(\tau)=\int_s^v\int_{|x|=\tau}|\nabla_xu|^2\,dxdt.
\]
\begin{lemma}\label{l:2.2}
Energy functions \eqref{2.7*} related to arbitrary solution
$u$ of problem \eqref{parabolicequation} satisfy the following
relationship:
\begin{multline}\label{2.16}
H(T,\tau)+I_s^{(T)}(\tau)\leq
c\,a(\tau)^{-\frac{2(1-\theta_2)}{2-(1-\theta_2)(1-q)}}
E(s,\tau)^{\frac2{2-(1-\theta_2)(1-q)}}\\+c_1\,
a(\tau)^{-\frac2{q+1}}E(s,\tau)^{\frac2{q+1}}+c\,a(\tau)^{-\frac2{q+1}}J_s^{(T)}(\tau)^{\frac2{q+1}}\\+
c\,a(\tau)^{-\frac{2(1-\theta_1)}{2-(1-\theta_1)(1-q)}}J_s^{(T)}(\tau)^{\frac2{2-(1-\theta_1)(1-q)}},\\
0<\theta_1=\tfrac{(q+1)+n(1-q)}{2(q+1)+n(1-q)}<1,\
\theta_2=\tfrac{n(1-q)}{2(q+1)+n(1-q)}.
\end{multline}
\end{lemma}
%\begin{proof}
Let us estimate the second term in right hand side of \eqref{2.1}.
By interpolation (see, for example,
\cite{DV85}) we have:
\begin{multline}\label{2.2}
\int_{|x|=\tau}|u|^2\, d\sigma\leq d_1\biggl(
\,\int_{\Omega(\tau)}|\nabla_xu|^2\,dx \biggr)^{\theta_1}\biggl(
\int _{\Omega(\tau)}|u|^{q+1}\,dx
\biggr)^{\frac{2(1-\theta_1)}{q+1}}\\+d_2\biggl(
\int_{\Omega(\tau)}|u|^{q+1}\,dx
\biggr)^{\frac{2}{q+1}}\quad\forall\,\tau>0,\ \theta_1\text{ is
from \eqref{2.16}}.
\end{multline}
Using \eqref{2.2} we easily arrive at
\begin{multline}\label{2.4}
\int_{|x|=\tau}|u||\nabla_xu|\,d\sigma \leq c\biggl(
\int_{|x|=\tau}|\nabla_xu|^2\,d\sigma \biggr)^{1/2}\left[ \biggl(
\int_{\Omega(\tau)}|\nabla_xu|^2\,dx
 \biggr)^{\frac{\theta_1}2}\right. \\
\left.\times \biggl( \int_{\Omega(\tau)}|u|^{q+1}
\biggr)^{\frac{1-\theta_1}{q+1}} + \biggl(
\int_{\Omega(\tau)}|u|^{q+1}\,dx \biggr)^{\frac1{q+1}} \right]=
c\biggl( \int_{|x|=\tau} |\nabla_xu|^2\,d\sigma
\biggr)^{1/2}\\\times\biggl( \int_{\Omega(\tau)}|\nabla_xu|^2\,dx
\biggr)^{\frac{\theta_1}2}\biggl( \int_{\Omega(\tau)}|u|^{q+1}\,dx
\biggr)^{\frac{1-\theta_1}{2}}\Biggl(
\int_{\Omega(\tau)}|u|^{q+1}\,dx
\biggr)^{\frac{(1-q)(1-\theta_1)}{2(q+1)}}\\+c\biggl(
\int_{|x|=\tau}|\nabla_xu|^2\,d\sigma \biggr)^{1/2}\biggl(
\int_{\Omega(\tau)}|u|^{q+1}\,dx \biggr)^{1/2}\biggl(
\int_{\Omega(\tau)}|u|^{q+1}\,dx \biggr)^{\frac{1-q}{2(q+1)}}.
\end{multline}
From condition \eqref{1.5} the monotonicity of
function $a(s)$ from \eqref{1.4} follows easily.
Therefore we can continue
estimating \eqref{2.4} as follows:
\begin{multline}\label{2.4*}
\int_{|x|=\tau}|u||\nabla_xu|\,d\sigma\leq c_1\biggl(
\int_{|x|=\tau}|\nabla_xu|^2\,d\sigma \biggr)^{1/2}\\\times\biggl(
\int_{\Omega(\tau)}|u|^2\,dx \biggr)^{\frac{(1-q)(1-\theta_1)}{4}}
a(\tau)^{-\frac{1-\theta_1}2}\biggl(
\int_{\Omega(\tau)}(|\nabla_xu|^2+a(|x|)|u^{q+1}|)\,dx
\biggr)^{1/2}\\+c_1a(\tau)^{-1/2}\biggl(
\int_{\Omega(\tau)}|u|^2\,dx
\biggr)^{\frac{1-q}4}\biggl(\int_{|x|=\tau}|\nabla_xu|^2\,d\sigma
\biggr)^{1/2}\biggl( \int_{\Omega(\tau)}a(|x|)|u|^{q+1}\,dx
\biggr)^{1/2}.
\end{multline}
Integrating \eqref{2.4*} in $t$ and using the Young inequality with
``$\varepsilon$'' we obtain:
\begin{multline}\label{2.5}
\int_s^{v}\int_{|x|=\tau}|u||\nabla u|\,d\sigma
dt\leq\varepsilon\int_{Q_s^{(v)}(\tau)}(|\nabla_xu|^2+a(|x|)|u|^{q+1})\,dxdt\\+c(\varepsilon)
a(\tau)^{-(1-\theta_1)}\sup_{s<t<v}\biggl(
\int_{\Omega(\tau)}|u(x,t)|^2\,dx
\biggr)^{\frac{(1-q)(1-\theta_1)}2}\int_s^{v}\int_{|x|=\tau}|\nabla_xu|^2\,d\sigma
dt\\+c(\varepsilon)\,a(\tau)^{-1}\sup_{s<t<v}\biggl(
\int_{\Omega(\tau)}|u(x,t)|^2\,dx
\biggr)^{\frac{1-q}2}\int_s^{v}\int_{|x|=\tau}|\nabla_xu|^2\,d\sigma
dt
\end{multline}
with arbitrary $v:s<v\leq T$. Let us fix now $v=\bar v=\bar
v(\tau,s)$ such that the following inequality holds:
\begin{equation}\label{2.6}
\int_{\Omega(\tau)}|u(x,\bar v)|^2\,dx\geq 2^{-1}\sup_{s\leq t\leq
T} \int_{\Omega(\tau)}|u(x,t)|^2\,dx
\end{equation}
Inserting inequality \eqref{2.5} with $v=\bar v$ into \eqref{2.1}
with $\hat t=\bar v$ and fixing ``$\varepsilon$'' small enough we
have:
\begin{multline}\label{2.7}
H(\bar v,\tau)+I_s^{(\bar v)}(\tau)\leq
H(s,\tau)+c\,a(\tau)^{-(1-\theta_1)}H(\bar
v,\tau)^{\frac{(1-q)(1-\theta_1)}2}J_s^{(\bar v)}(\tau)
\\+
c\,a(\tau)^{-1}\,H(\bar v,\tau)^{\frac{1-q}2}\,J_s^{(\bar
v)}(\tau),
\end{multline}
where $ J_s^{(v)}(\tau)$ is from \eqref{2.7*}. Using the Young
inequality again we deduce from \eqref{2.7}:
\begin{multline}\label{2.8}
H(\bar v,\tau)+I_s^{(\bar
v)}(\tau)\leq2H(s,\tau)+c\,a(\tau)^{-\frac{2(1-\theta_1)}{2-(1-\theta_1)(1-q)}}
(J_s^{(\bar v)}(\tau))^{\frac2{2-(1-\theta_1)(1-q)}}\\+
c\,a(\tau)^{-\frac2{1+q}}(J_s^{(\bar v)}(\tau))^{\frac2{1+q}}.
\end{multline}
Fixing now $v=T$ in \eqref{2.5} and using property \eqref{2.6} we
obtain the inequality:
\begin{multline}\label{2.9}
\int_s^T\int_{|x|=\tau}|u||\nabla_xu|\,d\sigma
dt\leq\varepsilon\,I_s^{(T)}(\tau)+c(\varepsilon)\,
a(\tau)^{-(1-\theta_1)}H(\bar v,\tau)^{\frac{(1-q)(1-\theta_1)}2}\,J_s^{(T)}(\tau)\\
+c(\varepsilon)\,a(\tau)^{-1}\,H(\bar
v,\tau)^{\frac{1-q}2}\,J_s^{(T)}(\tau).
\end{multline}
By $\hat{t}=T$ it follows from \eqref{2.1} due to \eqref{2.9} with
$\varepsilon=\frac12$:
\begin{multline}\label{2.10}
H(T,\tau)+I_s^{(T)}(\tau)\leq
H(s,\tau)+c\,a(\tau)^{-(1-\theta_1)}\,H(\bar
v,\tau)^{\frac{(1-q)(1-\theta)}2}\, J_s^{(T)}(\tau)\\*+
c\,a(\tau)^{-1}\,H(\bar v,\tau)^{\frac{1-q}2}\,J_s^{(T)}(\tau).
\end{multline}
From \eqref{2.8} we have
\begin{multline}
%\label{2.11}
H(\bar v,\tau)^\nu\leq
c\,H(s,\tau)^\nu+c\,a(\tau)^{-\frac{2(1-\theta_1)\nu}{2-(1-\theta_1)(1-q)}}
(J_s^{(\bar v)}(\tau))^{\frac{2\nu}{2-(1-\theta_1)(1-q)}}\\+
c\,a(\tau)^{-\frac{2\nu}{1+q}}(J_s^{(\bar
v)}(\tau))^{\frac{2\nu}{1+q}} \qquad \forall\, \nu>0.
\end{multline}
Using this estimate with $\nu_1=\frac{(1-q)(1-\theta_1)}2$ and
$\nu_2=\frac{1-q}2$  from \eqref{2.10} we deduce that
\begin{multline}\label{2.12}
H(T,\tau)+I_s^{(T)}(\tau)\leq
H(s,\tau)+c\,a(\tau)^{-(1-\theta_1)}\,H(s,\tau)^{\nu_1}\,J_s^{(T)}(\tau)\\+
c\,a(\tau)^{-1}\,H(s,\tau)^{\nu_2}\,J_s^{(T)}(\tau)+c\,a(\tau)^{-(1-\theta_1)\bigl(
1+\frac{2\nu_1}{2-(1-\theta_1)(1-q)}
\bigr)}\\\times(J_s^{(T)}(\tau))^{1+\frac{2\nu_1}{2-(1-\theta_1)(1-q)}}+
c\,a(\tau)^{-(1-\theta_1)-\frac{2\nu_1}{1+q}}\,
(J_s^{(T)}(\tau))^{1+\frac{2\nu_1}{1+q}}\\+
c\,a(\tau)^{-1-\frac{2(1-\theta_1)\nu_2}{2-(1-\theta_1)(1-q)}}\,(J_s^{(T)}(\tau))^{1+
\frac{2\nu_2}{2-(1-\theta_1)(1-q)}}\\+
c\,a(\tau)^{-1-\frac{2\nu_2}{1+q}}(J_s^{(T)}(\tau))^{1+\frac{2\nu_2}{1+q}}.
\end{multline}
Using the Young inequality we infer from \eqref{2.12}
\begin{multline}\label{2.13}
H(T,\tau)+I_s^{(T)}(\tau)\leq2H(s,\tau)+c\,a(\tau)^{-\frac2{1+q}}(J_s^{(T)}(\tau))^{\frac2{1+q}}
\\+c\,a(\tau)^{-\frac{2(1-\theta_1)}{2-(1-\theta_1)(1-q)}}(J_s^{(T)}(\tau))^{\frac2{2-(1-\theta_1)(1-q)}}.
\end{multline}
Now we have to estimate from above the term $H(s,\tau)$ in right
hand side of \eqref{2.13}. Due to the Gagliardo-Nirenberg
interpolation inequality  we have
\begin{multline}\label{2.14}
\int_{\Omega(\tau)}|u(x,s)|^2\,dx\leq d_3\biggl(
\int_{\Omega(\tau)}|\nabla_xu(x,s)|^2\,dx
\biggr)^{\theta_2}\biggl( \int_{\Omega(\tau)}|u(x,s)|^{q+1}
\biggr)^{\frac{2(1-\theta_2)}{q+1}}\\+ d_4\biggl(
\int_{\Omega(\tau)} |u(x,s)|^{q+1}\,dx \biggr)^{\frac2{q+1}},\quad
\theta_2\text{ is from \eqref{2.16}},
\end{multline}
and constants $d_3>0,\ d_4>0$ do not depend on $\tau$ as $\tau\to
0$. Taking into account the monotonicity of function $a(\tau)$ we
deduce from \eqref{2.14}
\begin{multline*}
\int_{\Omega(\tau)}|u(x,s)|^2\,dx\leq d_3\biggl(
\int_{\Omega(\tau)}|\nabla_xu(x,s)|^2dx
\biggr)^{\theta_2}\!\!\biggl(
\int_{\Omega(\tau)}a(|x|)|u(x,s)|^{q+1}dx \biggr)^{1-\theta_2}\\
\times
a(\tau)^{-(1-\theta_2)}\biggl(\int_{\Omega(\tau)}|u(x,s)|^{q+1}\,dx\biggr)^{\frac{(1-\theta_2)(1-q)}{1+q}}
\\+d_4\,a(\tau)^{-\frac2{q+1}}\biggl(
\int_{\Omega(\tau)}a(|x|)|u(x,s)|^{q+1}\,dx \biggr)^{\frac2{q+1}}
\\\leq
c\,a(\tau)^{-(1-\theta_2)}\int_{\Omega(\tau)}(|\nabla_xu|^2+a(|x|)|u(x,s)|^{q+1})\,dx\\*
\times\biggl( \int_{\Omega(\tau)}|u(x,s)|^2\,dx
\biggr)^{\frac{(1-\theta_2)(1-q)}2}+d_4\,a(\tau)^{-\frac2{q+1}}\biggl(
\int_{\Omega(\tau)}a(|x|)|u(x,s)|^{q+1}\,dx \biggr)^{\frac2{q+1}}.
\end{multline*}
Estimating the first term in the right hand side by the Young inequality with
``$\varepsilon$'', we have
\begin{multline}\label{2.15}
\int_{\Omega(\tau)}u^2(x,s)\,dx\leq
c_1\,a(\tau)^{-\frac2{q+1}}\Biggl(
\int_{\Omega(\tau)}a(|x|)|u(x,s)|^{q+1}\,dx
\Biggr)^{\frac2{q+1}}\\+
c\,a(\tau)^{-\frac{2(1-\theta_2)}{2-(1-\theta_2)(1-q)}}\Biggl(
\int_{\Omega(\tau)}(|\nabla_xu|^2+a(|x|)|u(x,s)|^{q+1})\,dx
\Biggr)^{\frac2{2-(1-\theta_2)(1-q)}}.
\end{multline}
Using \eqref{2.15} in \eqref{2.13} we obtain the
required \eqref{2.16}.
%\end{proof}

 Let us introduce the positive nondecreasing function
\begin{equation}\label{2.22}
s(\tau)=\tau^4\omega(\tau)^{-1},
\end{equation}
where $\omega(\tau)>0$ is from \eqref{1.1*}. Define the energy
function
\begin{equation}\label{y}
y(\tau)=I_{s(\tau)}^{(T)}(\tau),\quad\text{where }
I_{s}^{(T)}(\tau)\text{ is from }(\ref{2.7*}).
\end{equation}
\begin{lemma}\label{l:2.3}
The energy function $y(\tau)$ from \eqref{y} is the solution of the
following Cauchy problem for the ordinary differential inequality:
\begin{equation}\label{2.20}
y(\tau)\leq c_0\sum_{i=0}^2\Bigl( -\frac{y'(\tau)}{\psi_i(\tau)}
\Bigr)^{1+\lambda_i}\quad \forall\,\tau>0,
\end{equation}
\begin{equation}\label{2.21}
\quad y(0)\leq y_0, \quad y_0\text{ is from }(\ref{2.1*}),
\end{equation}
where
\begin{align*}
& \psi_0(\tau)=a(\tau)s'(\tau),\quad
\psi_1(\tau)=a(\tau)^{1-\theta_1},\quad
\psi_2(\tau)=a(\tau)^{1-\theta_2}s'(\tau),\\
&
\lambda_0=\frac{1-q}{1+q}>\lambda_2=\frac{(1-\theta_2)(1-q)}{2-(1-\theta_2)(1-q)}>\lambda_1=
\frac{(1-\theta_1)(1-q)}{2-(1-\theta_1)(1-q)}>0.
\end{align*}
\end{lemma}
%\begin{proof}
 It is easy to verify the following equality
\begin{multline}\label{2.17}
\frac d{d\tau}I_{s(\tau)}^{(T)}(\tau)=-\int_{s(\tau)}^T
\int_{\{|x|=\tau\}}(|\nabla_xu|^2+a(|x|)|u(x,s(\tau))|^{q+1})\,d\sigma
dt\\-
s'(\tau)\int_{\Omega(\tau)}(|\nabla_xu(x,s(\tau))|^2+a(|x|)|u(x,s(\tau))|^{q+1})\,dx
\end{multline}
Since $s'(\tau)\ge 0$, from \eqref{2.17} it follows that
\begin{equation}
%\label{2.18}
\int_{s(\tau)}^T\int_{\{|x|=\tau\}}|\nabla_xu|^2\,d\sigma\,dt=J_{s(\tau)}^{(T)}(\tau)\leq
-\frac d{d\tau}I_{s(\tau)}^{(T)}(\tau),
\end{equation}
\begin{equation}
%\label{2.19}
\int_{\Omega(\tau)}\!\!\bigl(
|\nabla_xu(x,s(\tau))|^2\!+\!a(|x|)|u(x,s(\tau))|^{q+1}
\bigr)dx\!=\!E(s(\tau),\tau)\!\leq\!-(s'(\tau))^{-1}\frac
d{d\tau}I_{s(\tau)}^{(T)}(\tau).
\end{equation}
Inserting these estimates in  \eqref{2.16} and
using additionally that $s'(\tau)\to 0$ as $\tau\to 0$
after simple calculations
we obtain
 ODI \eqref{2.20} and the initial condition
\eqref{2.21}.
%\end{proof}

Now we will study the asymptotic behavior of an arbitrary solution
$y(\tau)$ of system \eqref{2.20}, \eqref{2.21}. We have to prove
the existence of a continuous function $\bar\tau=\bar\tau(y_0)$ such
that $y(\tau)\leq0$ for arbitrary $\tau\geq\bar\tau(y_0)$.
Moreover, we have to find the sharp upper estimate for the function
$\bar\tau(y)$ as $y\to0$. It is
related to the optimal choice of the function $s(\tau)$, defined by
\eqref{2.22}. Consider the following auxiliary Cauchy problem:
\begin{equation}\label{2.23}
Y(\tau)=3c_0\max_{0\leq i\leq 2}\Bigl\{ \Bigl(
-\frac{Y'(\tau)}{\psi_i(\tau)} \Bigr)^{1+\lambda_i}\Bigr\},\quad
Y(0)=y_0>0,
\end{equation}
where $c_0>0$ is from \eqref{2.20}. It is easy to check the
following  comparison property:
\begin{equation}\label{2.24}
y(\tau)\leq Y(\tau)\quad \forall\, \tau>0,
\end{equation}
where $y(\tau)$ is arbitrary solution of the Cauchy problem
\eqref{2.20}, \eqref{2.21}.
\begin{lemma}\label{l:2.4}
Let $Y(\tau)$ be an arbitrary solution of the Cauchy problem
\eqref{2.23}. Then there exists a function $\bar\tau(r)<\infty \
\forall\,r>0$ such that $Y(\tau)\leq0\
\forall\,\tau>\bar\tau(y_0)$.
\end{lemma}

%\begin{proof}
Let us consider the following additional ordinary differential
equations (ODE):
\begin{equation}\label{2.25}
Y_i(\tau)=3c_0\Bigl(  -\frac{Y'_i(\tau)}{\psi_i(\tau)}
\Bigr)^{1+\lambda_i}, \quad i=0,1,2,
\end{equation}
or, equivalently:
\begin{equation}\label{2.25*}
Y'_i(\tau)=-\psi_i(\tau)\Bigl( \frac{Y_i(\tau)}{3c_0}
\Bigr)^{\frac1{1+\lambda_i}}:=-F_i(\tau,Y_i(\tau)).
\end{equation}
Let us define the following subdomains $\Omega_i\ i=0,1,2,$
\begin{align*}
& \Omega_0=\big\{  (\tau,y)\in\mathbb{R}_+^2:=\{ \tau>0,\ y>0
\}:F_0(\tau,y)=\min_{0\leq i\leq2}\{F_i(\tau,y)\} \big\},\\
& \Omega_1=\big\{  (\tau,y)\in\mathbb{R}_+^2 :
F_1(\tau,y)=\min_{0\leq i\leq2}\{F_i(\tau,y)\} \big\},\\
& \Omega_2=\big\{  (\tau,y)\in\mathbb{R}_+^2 :
F_2(\tau,y)=\min_{0\leq i\leq2}\{F_i(\tau,y)\} \big\}.
\end{align*}
It is easy to see that
$$
\Omega_0\cup\Omega_1\cup\Omega_2=\mathbb{R}_+^2.
$$
Due to \eqref{2.23}, \eqref{2.25}, \eqref{2.25*} it is easy to see
 that arbitrary solution $Y(\tau)$ of the problem \eqref{2.23} has the following structure:
\begin{equation}
\label{ad1} Y(\tau)=\big\{  Y_i(\tau)\quad\forall\,
(\tau,Y)\in\Omega_i,\ i=0,1,2 \big\},
\end{equation}
where $Y_i(\tau)$ is  solution of equation \eqref{2.25} (or
\eqref{2.25*}). It is easy to check that
\begin{align*}
& \Omega_0=\big\{  (\tau,y): y\geq 3c_0a(\tau)^{\frac2{1-q}} \big\},\\
& \Omega_1=\big\{  (\tau,y): y\leq 3c_0a(\tau)^{\frac2{1-q}}s'
(\tau)^{\frac{2}{(1-q)(\theta_1-\theta_2)}} \big\},\quad s'(\tau)
=\frac{ds(\tau)}{d\tau},\\
& \Omega_2=\big\{  (\tau,y):
3c_0a(\tau)^{\frac2{1-q}}s'(\tau)^{\frac{2}{(1-q)(\theta_1-\theta_2)}}\leq
y\leq 3c_0a(\tau)^{\frac2{1-q}} \big\}.
\end{align*}
Therefore  the solution $Y(\tau)$ of the Cauchy problem \eqref{2.23}
is dominated by the following curve:
\begin{equation}
\label{ad2} \tilde{Y}(\tau)=\begin{cases} y_0,&\quad \text{ if }\
0\leq\tau\leq\tau'\\
\tilde{Y}_2(\tau),&\quad\text{ if }\ \tau'\leq\tau\leq\tau''\\
\tilde{Y}_1(\tau),&\quad\text{ if }\ \tau''\leq\tau\leq\tau''',
\end{cases}
\end{equation}
where $\tau'$ is defined by equality
$y_0=3c_0a(\tau')^{\frac2{1-q}} \Rightarrow$
\begin{equation}
\label{ad3}
\frac{\tau'^2}{\omega(\tau')}=\frac2{1-q}(\ln(3c_0)-\ln y_0)^{-1},
\end{equation}
$\tilde{Y}_2(\tau)$ is the solution of the Cauchy problem:
\begin{equation}\label{2.29}
Y'_2(\tau)=-\psi_2(\tau)\Big(
\frac{Y_2(\tau)}{3c_0}\Big)^{\frac1{1+\lambda_2}},\quad
Y_2(\tau')=y_0,
\end{equation}
$\tau''$ is defined by the equality:
\begin{equation}\label{2.30}
\tilde{Y}_2(\tau'')=3c_0a(\tau'')^{\frac2{1-q}}s'(\tau'')^{\frac2{(1-q)(\theta_1-\theta_2)}}.
\end{equation}
Finally, $\tilde{Y}_1(\tau)$ is the solution of the Cauchy
problem:
\begin{equation}\label{2.31}
Y'_1(\tau)=-\psi_1(\tau)\Big(
\frac{Y_1(\tau)}{3c_0}\Big)^{\frac1{1+\lambda_1}},\quad
Y_1(\tau'')=\tilde{Y}_2(\tau''),
\end{equation}
and $\tau'''$ is such that $\tilde{Y}_1(\tau)\leq0\
\forall\,\tau\geq\tau'''$. It is easy to check that the solution of
\eqref{2.29} is
\begin{multline}
\label{ad4} \tilde Y_2(\tau)=\biggl[
y_0^{\frac{\lambda_2}{1+\lambda_2}}-\frac{\lambda_2}{(1+\lambda_2)(3c_0)^{\frac1{1+\lambda_0}}}
\int_{\tau'}^\tau  \psi_2(r)\,dr
\biggr]^{\frac{1+\lambda_2}{\lambda_2}}
\\*=\biggl[
y_0^{\frac{(1-\theta_2)(1-q)}2}-\frac{(1-\theta_2)(1-q)}{2(3c_0)^{\frac1{1+\lambda_0}}}
\int_{\tau'}^\tau a(r)^{1-\theta_2}s'(r)\,dr
\biggr]^{\frac2{(1-\theta_2)(1-q)}}.
\end{multline}
Equation \eqref{2.30} for $\tau''$ then yields:
\begin{multline}\label{2.33}
y_0^{\frac{(1-\theta_2)(1-q)}2}-\frac{(1-\theta_1)(1-q)}{2(3c_0)^{\frac1{1+\lambda_2}}}
\int_{\tau'}^{\tau''}a(r)^{1-\theta_2}s'(r)\,dr\\=(3c_0)^{\frac{(1-\theta_2)(1-q)}2}
a(\tau'')^{1-\theta_2}s'(\tau'')^2,\quad (\text{ since
}\tfrac{1-\theta_2}{\theta_1-\theta_2}=2).
\end{multline}
We will say that $a(\tau)\approx b(\tau)$, if there exist constant
$C$, which does not depend on $\tau$, such that
$$
0<C^{-1}a(\tau)\leq b(\tau)\leq Ca(\tau)\quad\forall\,
\tau:0<\tau<\tau_0.
$$
Due to condition \eqref{1.5} it follows easily too:
\begin{equation}\label{ss}
(2+\delta)\frac{\tau^3}{\omega(\tau)}\leq
s'(\tau)\leq\frac{4\tau^3}{\omega(\tau)}\quad\forall\,\tau>0.
\end{equation}
From definition \eqref{2.22} of $s(r)$ by virtue of \eqref{ss} and
Lemma A.1 we deduce
\begin{multline}\label{2.34}
\int_0^\tau a(r)^{1-\theta_2}s'(r)\,dr\approx\int_0^\tau\exp\Big(
-\frac{(1-\theta_2)\omega(r)}{r^2} \Big)r^3\omega(r)^{-1}\,dr\\
\approx\tau^6\omega(\tau)^{-2}\exp\Big(
-\frac{(1-\theta_2)\omega(\tau)}{\tau^2} \Big)\approx
a(\tau)^{1-\theta_2}(s'(\tau))^2\quad \forall\,
\tau:0<\tau<\tau_0<\infty.
\end{multline}
Thus, from \eqref{2.33} due to \eqref{2.34} one obtains the following estimate for
$\tau''$
\begin{equation}\label{2.35}
c_1y_0^{\frac{(1-\theta_2)(1-q)}2}\leq
a(\tau'')^{1-\theta_2}s'(\tau'')^2\leq
c_2y_0^{\frac{(1-\theta_2)(1-q)}2},
\end{equation}
where positive constants $c_1,\ c_2$ does not depend on $y_0$.
Now, the solution of the Cauchy problem \eqref{2.31} is:
\begin{equation}
%\label{2.36}
\tilde Y_1(\tau)=\biggl[
\tilde{Y}_2(\tau'')^{\frac{(1-\theta_1)(1-q)}2}-\frac{(1-\theta_1)(1-q)}{2(3c_0)^{\frac1{1+\lambda_1}}}
\int _{\tau''}^\tau a(r)^{1-\theta_1}\,dr
\biggr]^{\frac2{(1-\theta_1)(1-q)}}.
\end{equation}
Thus, $\tau'''$ is defined by the equation:
\begin{equation}\label{2.37}
\tilde{Y}_2(\tau'')^{\frac{(1-\theta_1)(1-q)}2}-\frac{(1-\theta_1)(1-q)}{2(3c_0)^{\frac1{1+\lambda_1}}}
\int_{\tau''}^{\tau'''}a(r)^{1-\theta_1}\,dr.
\end{equation}
Due to Lemma A.1 we have
\begin{multline}\label{2.38}
\biggl( \int_0^\tau a(r)^{1-\theta_1}\,dr \biggr)^2\approx\biggl(
 \int_0^\tau \exp\Big( -\frac{\beta\omega(r)}{r^2} \Big)dr
 \biggr)^2\\\approx \biggl( \frac{\tau^3}{\omega(\tau)}\exp\Bigl( -\frac{\beta\omega(\tau)}{\tau^2}
  \Bigr) \biggr)^2\approx s'(\tau)^2a(\tau)^{1-\theta_2}\quad
  \forall\, \tau>0,
\end{multline}
where
$1-\theta_1=\beta=\frac{q+1}{2(q+1)+u(1-q)}=\frac{1-\theta_2}2$.
It is easy to see that
$$
\int_0^\tau a(r)^{1-\theta_1}\,dr\approx\int_{\frac\tau2}^\tau
a(r)^{1-\theta_1}\,dr\quad \text{if}\quad \tau\to 0.
$$
Therefore due to \eqref{2.37} the following inequalities are
sufficient conditions for  $\tau'''$:
$$
a(\tau''')^{1-\theta_2}s'(\tau''')^2\leq c_3\tilde
Y_2(\tau'')^{\frac{(1-\theta_2)(1-q)}2},\quad \tau'''>2\tau''.
$$
Finally, by virtue of \eqref{ad4} we obtain the following unique
sufficient condition which defines $\tau'''$:
\begin{equation}\label{2.39}
a(\tau''')^{1-\theta_2}s'(\tau''')^2\leq
c_4y_0^{\frac{(1-\theta_2)(1-q)}2},
\end{equation}
Condition \eqref{2.39} can be rewritten in the form:
\begin{equation}\label{2.40}
\frac{\exp\Bigl(
 -\frac{(1-\theta_2)(1-\nu)\omega(\tau''')}{(\tau''')^2}\Bigr)
 \cdot\exp\Bigl( -(1-\theta_2)\nu\frac{\omega(\tau''')}{(\tau''')^2} \Bigr)
  \omega(\tau''')}
 {\big(\frac{\omega(\tau''')}{(\tau''')^2}\big)^3}\leq c_5y_0^{\frac{(1-\theta_2)(1-q)}2}
\end{equation}
with arbitrary $1>\nu>0$. It is obviously, that the following is a sufficient
condition for \eqref{2.40}
$$
\exp\Bigl( -\frac{(1-\theta_2)(1-\nu)\omega(\tau''')}{(\tau''')^2}
\Bigr)\leq c_6y_0^{\frac{(1-\theta_2)(1-q)}2},\quad
c_6=c_6(\nu,\omega_0,c_5)
$$
or,
\begin{equation}
%\label{2.41}
\frac{(\tau''')^2}{\omega(\tau''')}\leq c_7(\ln
y_0^{-1})^{-1},\quad c_7=c_7(c_6,\nu,\omega_0), \quad \omega_{0}
\quad \text{is from} \quad (A_3).
\end{equation}
Thus, the assertion of Lemma \ref{l:2.4} holds with
$\bar\tau(r)$ defined by:
\begin{equation}\label{ad5}
\frac{\bar\tau(r)^2}{\omega(\bar\tau(r))}=c_7(\ln
r^{-1})^{-1}\quad\forall\,r>0.
\end{equation}
%\end{proof}

%\begin{proof}[Proof of Theorem 1.1.]
\noindent
{\it Proof of Theorem 1.1.}
Due to Lemma A.3 from Appendix we can suppose that
\begin{equation}
%\label{2.43}
y_0\ll1\quad\text{and}\quad \bar\tau(y_0)<1.
\end{equation}
From definition \eqref{2.20} of function $y(\tau)$ due to Lemma
\ref{l:2.4} and property \eqref{2.24} it follows that
$$
I_{s(\bar\tau(y_0))}^{(T)}(\bar\tau(y_0))=0\quad\text{for
arbitrary}\quad T<\infty.
$$
Therefore our solution $u(x,t)$ has the following property:
\begin{equation}\label{2.44}
u(x,t)\equiv0 \quad\forall\, (x,t)\in\big\{ |x|\geq\tau_1,\ t\geq
s(\tau_1)\big\},\quad\tau_1=\bar\tau(y_0).
\end{equation}
From identity \eqref{2.1} with $\tau=0$ we deduce that
\begin{equation}\label{2.45}
\frac
d{dt}\int_{\Omega}|u(x,t)|^2\,dx+\int_{\Omega}(|\nabla_xu(x,t)|^2+a_0(x)|u|^{q+1})\,dx\leq0\quad
\forall\, t\in(s(\tau_1),T).
\end{equation}
Due to \eqref{2.44} and the Poincar{\'{e} inequality it follows from
\eqref{2.45}:
\begin{equation}\label{2.46}
H'(t)+\frac{\bar
c}{\tau_1^2}H(t)\leq0\quad\forall\,t>s(\tau_1),\quad \bar
c=\text{const}>0,
\end{equation}
where $H(t):=H(t,0)$, $H(t,\tau)$ is defined by \eqref{2.4*},
constant $\bar c>0$ does not depend on $t$. Integrating ODI
\eqref{2.46} we deduce the following relationship easily:
$$
H(t+s(\tau_1))\leq H(s(\tau_1))\exp\Bigl( -\frac{\bar
ct}{\tau_1^2} \Bigr)\quad\forall\, t>0.
$$
Using additionally estimate \eqref{2.1*} with $\hat{t}=s(\tau_1)$
we deduce:
\begin{equation}\label{2.47}
H(t+s(\tau_1))\leq  y_0\exp\Bigl( -\frac{\bar ct}{\tau_1^2}
\Bigr)\quad\forall\,t>0.
\end{equation}
Define  $t_1>0$ by
\begin{equation}\label{2.48}
y_0\exp\Bigl( -\frac{\bar ct_1}{\tau_1^2}
\Bigr)=y_0^{1+\gamma}\quad\Leftrightarrow\quad t_1=\frac{\gamma\ln
y_0^{-1}}{\bar c}\tau_1^2,\quad \gamma=\text{const}>0.
\end{equation}
Due to \eqref{ad5} from \eqref{2.48}  it follows that
\begin{equation}\label{2.49}
t_1=\frac{\gamma c_7}{\bar c}\omega(\tau_1).
\end{equation}
Thus, we have:
\begin{equation}\label{2.50}
H(t_1+s(\tau_1))=\int_{\Omega}|u(x,t_1+s(\tau_1))|^2\,dx\leq
y_0^{1+\gamma},\quad \gamma>0.
\end{equation}
So, we finished first round of computations. For the second round we
will consider our initial-boundary problem
\eqref{parabolicequation} in the domain
$\Omega\times(t_1+s(\tau_1),\infty)$ with initial data
\eqref{2.50} instead of \eqref{2.1*}. Repeating all previous
computations we deduce the following analogue of estimate \eqref{2.50}
\begin{equation}\label{2.51}
H(t_2+s(\tau_2)+t_1+s(\tau_1))\leq y_0^{(1+\gamma)^2},
\end{equation}
where  as in \eqref{ad5} and \eqref{2.48}
\begin{equation}
%\label{2.52}
\tau_2^2=c_7\omega(\tau_2)(\ln
y_0^{-(1+\gamma)})^{-1}=\frac{c_7}{1+\gamma}\omega(\tau_2)(\ln
y_0^{-1})^{-1},\quad \tau_2=\bar\tau(y_0^{1+\gamma}).
\end{equation}
Analogously to \eqref{2.49} we have also:
\begin{equation}
%\label{2.53}
t_2=\frac{\gamma\ln y_0^{-(1+\gamma)}}{\bar
c}\tau_2^2=\frac{\gamma c_7}{\bar c}\omega(\tau_2).
\end{equation}
Now using estimate \eqref{2.51} as a starting point for next round
of computations we find $\tau_3,\ t_3$ and so on. As result, after
$j$ rounds we get
\begin{equation}\label{2.54}
H\biggl(\,  \sum_{i=1}^j t_i+\sum_{i=1}^j s(\tau_i) \biggr)\leq
y_0^{(1+\gamma)^j}\to 0 \quad \text{as}\quad j\to 0,
\end{equation}
where
\begin{equation}\label{2.55}
\tau_i^2\leq\frac{c_7\omega(\tau_i)}{(1+\gamma)^{i-1}}(\ln
y_0^{-1})^{-1}
\end{equation}
Due to condition $(A_3)$ it follows from \eqref{2.55}:
\begin{equation*}
\tau_i^2\leq\frac{c_7\omega_0 (\ln
y_0^{-1})^{-1}}{(1+\gamma)^{i-1}}
\end{equation*}
From definition \eqref{2.22} of function $s(\tau)$ due to
condition \eqref{1.5} it follows the estimate
$$
s(\tau)\leq\tau_0^2\omega(\tau_0)^{-1}\tau^2\quad\forall\,\tau_0>0,\
\forall\,\tau>0.
$$
Therefore inequality \eqref{2.55} yields:
\begin{equation}\label{2.56}
\sum_{i=1}^\infty s(\tau_i)<\tilde{c}<\infty.
\end{equation}
Obviously, we have also: $t_i=\frac{\gamma c_7}{\bar
c}\omega(\tau_i)$. Therefore, due to \eqref{2.55} we have:
\begin{equation}
%\label{2.57}
\sum_{i=1}^j t_i=\frac{\gamma c_7}{\bar c}\sum_{i=1}^j
\omega(\tau_i)\leq C\sum_{i=1}^j\omega(C_1\lambda^i),
\end{equation}
where $C=\frac{\gamma c_7}{\bar c}, \ C_1=\Bigl(
\frac{c_7\omega_0}{\ln y_0^{-1}(1+\gamma)} \Bigr)^{1/2},\
\lambda=(1+\gamma)^{-1/2}<1$. In virtue of condition \eqref{1.2}
it is easy to check that
\begin{equation}\label{2.58}
\sum_{i=1}^j \omega(C_1\lambda^i)\approx
\ln\lambda^{-1}\int_{C_1\lambda^j}^{C_1}\frac{\omega(s)}s\,ds<c<\infty\quad
\forall\,j\in\mathbb N.
\end{equation}
From \eqref{2.54} due to \eqref{2.56}, \eqref{2.58} and condition
\eqref{1.2} it follows that
$$
H(R)=0, \quad R=\sum_{i=1}^\infty t_i+\sum_{i=1}^\infty
s(\tau_i)<\infty,
$$
which completes the proof of Theorem~\ref{MR}. \qed
%\end{proof}
%
\section{\bf Dini condition \eqref{1.2} of extinction in finite
time via semi-classical limit of Schr\"odinger operator}
\def\theequation{3.\arabic{equation}}\makeatother
\setcounter{equation}{0}

Here we prove Proposition 1.3. We recall the definition of
$\lambda_1(h)$ and $\mu(\alpha)$ for $h>0$ and $\alpha>0$ :
\[
\lambda_1(h) = \inf \left\{ \int_{B_1} |\nabla v|^2 + h^{-2} a(|x|) |v|^2 \;
dx : \; v \in W^{1,2}(B_1), \; ||v||_{L_2(B_1)}=1 \right\},
\]
and
\[
\mu(\alpha) = \lambda_1(\alpha^\frac{1-q}{2}).
\]
We define $r(z)=a^{-1}(z)$ or equivalently $z=a(r(z))$ and $\rho(z)=z(r(z))^2$ for $z$ small enough.
 We will use the following technical statement
\begin{lemma}[\rm Corollaries 2.23, 2.31 in \cite{Be04}]
Under assumptions $(A_1)-(A_3)$\\ and $\eqref{1.7}$, there exist
four positives constants $C_1$, $C_2$, $C_3$ and $C_4$ such that
for $h$ small enough,
\[
C_1 h^{-2} \rho^{-1}(C_2 h^2) \leq \lambda_1(h) \leq C_3 h^{-2}
\rho^{-1}(C_4 h^2).
\]
\end{lemma}
Our main starting point in the proof of Proposition 1.3 is the
following

\medskip
\noindent {\bf Theorem A} (Th. 2.2 in \cite{BHV01}). {\it Under
assumptions $(A_1)-(A_3)$, if there exists a decreasing sequence
$(\alpha_n)$  of positive real numbers such that
\[
\sum_{n=0}^{+\infty} \frac{1}{\mu(\alpha_n)} \left(
\ln(\mu(\alpha_n))+\ln\left(\frac{\alpha_n}{\alpha_{n+1}}\right)+1
\right) < +\infty,
\]
then problem (\ref{parabolicequation}) satisfies the TCS-property.}

\medskip

The first step in the proof of Proposition 1.3 is the estimation
of $\rho^{-1}$ in a neighbourhood of zero.

\begin{lemma}
Under assumptions $(A_1)-(A_3)$ with $\eqref{1.6}$ there holds
\begin{equation} \label{cestimate}
\frac{s}{(1+\alpha)} \ln \left(\frac{1}{s}\right)
\frac{1}{\omega\left( \left( \frac{\omega_0(1+\alpha)}{\ln
\left(\frac{1}{s}\right)}\right)^\frac{1}{2} \right)} \leq
\rho^{-1}(s) \leq s \ln \left(\frac{1}{s}\right)
\frac{1}{\omega\left( \left( \frac{1}{\ln
\left(\frac{1}{s}\right)} \right)^\frac{1}{\delta} \right)},
\end{equation}
for arbitrary $\alpha >0$, for all $s>0$ small enough.
\end{lemma}
%\begin{proof}
First of all, we prove the following estimate for $\rho(z)$:
\begin{eqnarray} \label{estimateforrho}
\!\!\!\!\!\! \rho(z)\ln \left(\frac{1}{z}\right)
\left[ \omega\left(\left(\frac{\omega_0}{\ln \left(\frac{1}{z}\right)}\right)^\frac{1}{2}\right)\right]^{-1}
  \leq z \leq \rho(z)\ln \left(\frac{1}{z}\right)
  \left[{\omega\Bigg( \bigg( \frac{1}{\ln \left(\frac{1}{z}\right)} \bigg)^\frac{1}{\delta}
  \Bigg)}\right]^{-1}.
\end{eqnarray}
Starting with $r >0$ small enough, we have from \eqref{1.6} the relationship\\ $r^{2-\delta} \leq \omega(r)
\leq \omega_0$ and since for $z>0$ small enough,
\[
(r(z))^2 \ln \left(\frac{1}{z}\right)=\omega(r(z))\Longrightarrow
r(z)^{2-\delta} \leq (r(z))^2 \ln \left(\frac{1}{z}\right) \leq
\omega_0.
\]
Therefore, we obtain
\begin{eqnarray} \label{estimateforr}
\left( \frac{1}{\ln \left(\frac{1}{z}\right)}
\right)^\frac{1}{\delta} \leq r(z) \leq \left(\frac{\omega_0}{\ln
\left(\frac{1}{z}\right)}\right)^\frac{1}{2}.
\end{eqnarray}
Since $\omega$ is a non decreasing function,
\[
\omega\left( \left( \frac{1}{\ln \left(\frac{1}{z}\right)}
\right)^\frac{1}{\delta} \right) \leq \omega(r(z)) \leq
\omega\left(\left(\frac{\omega_0}{\ln
\left(\frac{1}{z}\right)}\right)^\frac{1}{2}\right).
\]
Substituting the definition of $\omega(r)$,
\[
\omega\left( \left( \frac{1}{\ln \left(\frac{1}{z}\right)}
\right)^\frac{1}{\delta} \right) \leq (r(z))^2 \ln
\left(\frac{1}{z}\right) \leq
\omega\left(\left(\frac{\omega_0}{\ln
\left(\frac{1}{z}\right)}\right)^\frac{1}{2}\right).
\]
It follows the estimate for $\rho(z)$.
\begin{eqnarray} \label{estimateforrhoz}
z \frac{1}{\ln \left(\frac{1}{z}\right)} \omega\left( \left( \frac{1}{\ln \left(\frac{1}{z}\right)}
 \right)^\frac{1}{\delta} \right) \leq \rho(z) \leq z \frac{1}{\ln \left(\frac{1}{z}\right)} \omega\left(\left(\frac{\omega_0}{\ln \left(\frac{1}{z}\right)}\right)^\frac{1}{2}\right).
\end{eqnarray}
By an easy calculation, we have (\ref{estimateforrho}).\\
\vspace{0 cm}\\
Here and further, $z=\rho^{-1}(s)$. By using
\eqref{estimateforr} and $\rho(z)=z(r(z))^2$,
\[
\rho(z) \geq z \left(\frac{1}{\ln
\left(\frac{1}{z}\right)}\right)^\frac{2}{\delta}
\Longleftrightarrow \frac{1}{\rho(z)} \leq \frac{1}{z} \left( \ln
\left(\frac{1}{z}\right) \right)^\frac{2}{\delta},
\]
or equivalently,
\[
\ln \left(\frac{1}{\rho(z)}\right) \leq \ln
\left(\frac{1}{z}\right) + \frac{2}{\delta} \ln \left( \ln
\left(\frac{1}{z}\right) \right).
\]
Let $\alpha>0$. Then for $z$ small enough, since $\ln(\ln(z^{-1})) << \ln z^{-1}$,
\begin{eqnarray} \label{estimatefrombelowforrhominusone}
\ln \left(\frac{1}{\rho(z)}\right) \leq (1+\alpha)\ln
\left(\frac{1}{z}\right)\Longleftrightarrow \rho(z) \geq
z^{1+\alpha}\Longrightarrow \rho^{-1}(s) \leq
s^\frac{1}{1+\alpha}.
\end{eqnarray}
Substituting $z=\rho^{-1}(s)$ in \eqref{estimateforrho} yields
\[
s \ln \left(\frac{1}{\rho^{-1}(s)}\right)
\left[{\omega\Bigg(\bigg(\frac{\omega_0}{\ln
\left(\frac{1}{\rho^{-1}(s)}\right)}\bigg)^\frac{1}{2}\Bigg)}\right]^{-1}
\leq \rho^{-1}(s),
\]
and due to \eqref{estimatefrombelowforrhominusone},
\[
\frac{s}{(1+\alpha)} \ln \left(\frac{1}{s}\right)
\left[{\omega\Bigg(\bigg(\frac{\omega_0(1+\alpha)}{\ln
\left(\frac{1}{s}\right)}\bigg)^\frac{1}{2}\Bigg)}\right]^{-1}
\leq \rho^{-1}(s),
\]
since $\omega$ is a nondecreasing function.\\
For the right-hand side of \eqref{cestimate}, we substitute $z=\rho^{-1}(s)$ in \eqref{estimateforrho}.
\[
\rho^{-1}(s) \leq s \ln \left(\frac{1}{\rho^{-1}(s)}\right)
\left[{\omega\Bigg( \bigg( \frac{1}{\ln
\left(\frac{1}{\rho^{-1}(s)}\right)} \bigg)^\frac{1}{\delta}
\Bigg)}\right]^{-1}.
\]
But from (\ref{estimateforr}), $r(z) \to z$ so we have for $z$ small enough,
$\rho(z) \leq z,$
which gives
$\rho^{-1}(s) \geq s.$
Consequently,
\[
\rho^{-1}(s) \leq s \ln \left(\frac{1}{s}\right)
\left[{\omega\Bigg( \bigg( \frac{1}{\ln \left(\frac{1}{s}\right)}
\bigg)^\frac{1}{\delta} \Bigg)}\right]^{-1},
\]
which completes the proof.
%\end{proof}
\begin{lemma}
Under $(A_1)-(A_3)$ with $\eqref{1.6}$ and $\eqref{1.7}$, if
\begin{equation}\label{series}
\sum_{n=n_0}^{+\infty} \frac{\omega\left(\frac{1}{\left(n\ln n\right)^\frac{1}{2}}\right)
}{n} < + \infty,
\end{equation}
then all solutions of $(\ref{parabolicequation})$ vanish in a
finite time. Moreover,
\begin{equation}\label{seriesintegral}
\sum_{n=n_0}^{+\infty} \frac{\omega\left(\frac{1}{\left(n\ln n\right)^\frac{1}{2}}\right)
}{n} < + \infty \; \Longleftrightarrow \; \int_0^c \frac{\omega(x)}{x} \, dx < + \infty.
\end{equation}
\end{lemma}
%\begin{proof}
From Lemma 3.3 and Lemma 3.1 we get
\[
K_1 \ln \left(\frac{1}{h}\right) {\omega\left(\frac{K_2}{\left(\ln
\left(\frac{1}{h}\right)\right)^\frac{1}{2}}\right)} \leq
\lambda_1(h) \leq K_3 \ln \left(\frac{1}{h}\right) \left
[{\omega\left( \frac{K_4}{\left( \ln \left(\frac{1}{h}\right)
\right)^\frac{1}{\delta} } \right)}\right ]^{-1},
\]
and since $\omega(r) \geq r^\theta$ for $r$ small enough, we have
\[
K_1 \ln \left(\frac{1}{h}\right)
\left[{\omega\Bigg(\frac{K_2}{\left(\ln
\left(\frac{1}{h}\right)\right)^\frac{1}{2}}\Bigg)}\right]^{-1}
\leq \lambda_1(h) \leq K'_3 \ln
\left(\frac{1}{h}\right)^{1+\frac{2-\delta}{\delta}},
\]
which leads to
\begin{eqnarray}
C'_1 \ln \left(\frac{1}{h}\right) \left
[{\omega\left(\frac{C'_2}{\left(\ln
\left(\frac{1}{h}\right)\right)^\frac{1}{2}}\right)}\right ]^{-1}
\leq \lambda_1(h) \leq C'_3 \ln
\left(\frac{1}{h}\right)^\frac{2}{\delta}.
\end{eqnarray}
The real number $\alpha$ is defined by $h=\alpha^\frac{1-q}{2}$ and thus,
\[
C''_1 \ln \left(\frac{1}{\alpha}\right)
\left[{\omega\Bigg(\frac{C''_2}{\left(\ln\left(\frac{1}{\alpha}\right)\right)^\frac{1}{2}}\Bigg)}\right]^{-1}
\leq \mu(\alpha) \leq C''_3 \ln
\left(\frac{1}{\alpha}\right)^\frac{2}{\delta}.
\]
From Theorem A, if $(\alpha_n)$ is a decreasing sequence of positive real numbers and
\[
\sum_{n=n_0}^{+\infty} \frac{\omega\left(\frac{C''_2}{\left(\ln\left(\frac{1}{\alpha_n}\right)\right)^\frac{1}{2}}\right)
}{\ln \left(\frac{1}{\alpha_n}\right)} \left[ \ln\left(\ln\left(\frac{1}{\alpha_n}\right)\right)
+ \ln \left(\frac{\alpha_n}{\alpha_{n+1}} \right) + 1 \right] < + \infty,
\]
then all the solutions of $(\ref{parabolicequation})$ vanish in a finite time.\\
The main point is the sequence $(\alpha_n)$. In \cite{BHV01}, they set $\alpha_n=2^{-n}$.
A better choice is $\alpha_n=n^{-Kn}$ for some $K > 0$ since
$\displaystyle \ln\left(\ln\left(\frac{1}{\alpha_n}\right)\right) \sim \ln
\left(\frac{\alpha_n}{\alpha_{n+1}} \right)$ which leads to \eqref{series}.\\
Now, we have to show that
\[
\sum_{n=n_0}^{+\infty} \frac{\omega\left(\frac{1}{\left(n\ln n\right)^\frac{1}{2}}\right)
}{n} < + \infty \; \Longleftrightarrow \; \int_0^c \frac{\omega(x)}{x} \, dx < + \infty.
\]
The series is finite if and only if
$\displaystyle \int_{n_0}^{+ \infty}
\frac{\omega\left(\frac{1}{\left(x\ln x\right)^\frac{1}{2}}\right)}{x} \, dx =
\int_0^{1/n_0} \frac{\omega\left(\left(\frac{x}{-\ln x}\right)^\frac{1}{2}\right)}{x} \, dx$
is finite. The following inequalities hold for $c>0$ small enough:
\[
\int_0^c \frac{\omega(x)}{x} \, dx \leq \int_0^c
\frac{\omega\left(\left(\frac{x}{-\ln x}\right)^\frac{1}{2}\right)}{x} \, dx
\leq \int_0^c \frac{\omega(x^\frac{1}{2})}{x} \, dx = 2 \int_0^{\sqrt{c}} \frac{\omega(x)}{x} \, dx,
\]
which completes the proof of Proposition 1.3.
%\end{proof}
\section{Appendix}
\renewcommand{\theequation}{4.\arabic{equation}}
\setcounter{section}{1}
\renewcommand{\thesection}{\Alph{section}}
\setcounter{equation}{0}
\begin{lemma}
Let the nonnegative nondecreasing function $\omega(s),\ s\geq0$,
satisfy  condition \eqref{1.5}. Then for any
$m\in\mathbb{R}^1,\ l\in\mathbb{R}^1,\ A>0$, one has
\begin{equation}\label{4.1}
\int_0^\tau s^{m-2}\omega(s)^{l+1}\exp\Bigl(
-\frac{A\omega(s)}{s^2}
\Bigr)ds\approx{\tau^{m+1}\omega(\tau)^l}\exp\Bigl(
-\frac{A\omega(\tau)}{\tau^2} \Bigr)\quad\text{as}\quad \tau\to0.
\end{equation}
\end{lemma}
%\begin{proof}
It is easy to check the following equality
\begin{multline}\label{4.1*}
\frac{d}{ds}\Big(s^{m+1}\omega(s)^l\exp\Big(-\frac{A\omega(s)}{s^2}\Big)\Big)
=s^m\omega(s)^l\exp\Big(-\frac{A\omega(s)}{s^2}\Big)\\\times \Big[
(m+1)+l\frac{s\omega'(s)}{\omega(s)}+\frac{A\omega(s)}{s^2}\Big(2-\frac{s\omega'(s)}{\omega(s)}\Big)\Big]
\equiv
s^m\omega(s)^l\exp\Big(-\frac{A\omega(s)}{s^2}\Big)[I_1+I_2+I_3].
\end{multline}
Integrating condition \eqref{1.5} we get:
\begin{equation}\label{4.1**}
\omega(s)\geq s^{2-\delta}\quad\forall\,s\in(0,s_0),
\end{equation}
and, as a consequence, $\frac{\omega(s)}{s^2}\to\infty$ as $s\to 0$.
Now due to \eqref{1.5} it follows that
$$
I_3\gg|I_1|,\quad I_3\gg|I_2|\quad \text{as }s\to 0.
$$
Therefore, integrating \eqref{4.1*} we obtain \eqref{4.1}.
%\end{proof}
\begin{lemma}
Let $\Omega$ be a domain from problem \eqref{parabolicequation}, let
$\Omega_0$ be a subdomain of
$\Omega:\overline{\Omega}_0\subset\Omega$. Then the following
interpolation inequality holds
\begin{equation}\label{4.2}
\biggl( \int_\Omega v^2(x)\,dx \biggr)^{1/2} \leq c_1\biggl(
\int_\Omega |\nabla_xv|^2\,dx \biggr)^{1/2}+c_2\biggl(
\int_{\Omega_0} |v|^\lambda\,dx \biggr)^{1/\lambda}\quad \forall\,
v\in W_2^1(\Omega),
\end{equation}
where $\lambda:1<\lambda\leq2$, positive constants $c_1,\ c_2$
does not depend on $v$.
\end{lemma}
%\begin{proof}
We start from the standard interpolation inequality
\begin{equation}\label{4.3}
\biggl( \int_\Omega v^2(x)\,dx \biggr)^{1/2} < c_1\biggl(
\int_\Omega |\nabla_xv|^2\,dx \biggr)^{1/2}+c_2\biggl(
\int_{\Omega} |v|^\lambda\,dx \biggr)^{1/\lambda}\quad \forall\,
v\in W_2^1(\Omega),
\end{equation}
It is clear that
\begin{equation}\label{4.4}
\biggl(  \int_\Omega |v|^\lambda\,dx
\biggr)^{1/\lambda}\leq\biggl(  \int_{\Omega_0}|v|^\lambda\,dx
\biggr)^{1/\lambda}+\biggl(
\int_{\Omega\setminus\Omega_0}|v|^\lambda\,dx \biggr)^{1/\lambda}.
\end{equation}
Let $\Omega'_0$ be a subdomain of $\Omega_0$ such that
$\overline{\Omega'}_0\subset\Omega_0$. let $\xi(x)\geq0$ be
$C^1$-smooth function such that
\begin{equation}\label{4.5}
\xi(x)=0\quad\forall\,x\in\overline{\Omega'}_0,\qquad
\xi(x)=1\quad\forall\,x\in\Omega\setminus\Omega_0.
\end{equation}
Then we have due to the Poincar\'{e} inequality:
\begin{multline}\label{4.6}
\biggl( \int_{\Omega\setminus\Omega_0}|v|^\lambda\,dx
\biggr)^{1/\lambda}\leq\biggl(
\int_{\Omega\setminus\Omega'_0}|v\xi|^\lambda\,de
\biggr)^{1/\lambda}\\\leq c\biggl(
\int_{\Omega\setminus\Omega'_0}|\nabla_x(b\xi)|^\lambda\,dx
\biggr)^{1/\lambda}\leq c\biggl(
\int_{\Omega\setminus\Omega'_0}|\nabla v |^\lambda\,dx
\biggr)^{1/\lambda}\\+c\biggl(
\int_{\Omega_0\setminus\Omega'_0}|\nabla\xi|^\lambda|v|^\lambda\,dx
\biggr)^{1/\lambda}\leq c_1\biggl(
\int_{\Omega\setminus\Omega'_0}|\nabla v|^2\,dx
\biggr)^{1/2}+c_2\biggl(
\int_{\Omega_0\setminus\Omega'_0}|v|^{\lambda}\,dx
\biggr)^{1/\lambda}.
\end{multline}
From \eqref{4.3} due to \eqref{4.4}--\eqref{4.6} one obtains
\eqref{4.2}. Lemma A.2 is proved.
%\end{proof}
\begin{lemma}
Let $u(x,t)$ be an arbitrary energy solution of problem
\eqref{parabolicequation}. Then
$H(t)=\int_{\Omega}|u(x,t)|^2\,dx\to0$ as $t\to\infty$.
\end{lemma}
%\begin{proof}
It is clear that there exists a constant $a_0>0$ and a subdomain
$\Omega_0\subset\Omega$ such that $a(x)\geq a_0>0$ for all
$x\in\overline{\Omega}_0$. From \eqref{2.47} it follows that
\begin{equation}\label{4.7}
\frac
d{dt}\int_\Omega|u(x,t)|^2\,dx+\int_{\Omega}|\nabla_xu(x,t)|^2\,dx
+a_0\int_{\Omega_0}|u(x,t)|^{q+1}\,dx\leq0.
\end{equation}
Due to Lemma A.2 we have
\begin{equation}\label{4.8}
\varepsilon\int_{\Omega}|u|^2\,dx\leq\varepsilon
c_1\int_{\Omega}|\nabla_xu|^2\,dx+\varepsilon c_2\biggl(
\int_{\Omega_0}|u|^{q+1}\,dx \biggr)^{\frac2{1+q}}\quad
\forall\,\varepsilon>0.
\end{equation}
Adding \eqref{4.7} and \eqref{4.8} we get
\begin{multline}\label{4.9}
\frac{d}{dt}\int_{\Omega}|u(x,t)|^2\,dx+\varepsilon\int_{\Omega}|u(x,t)|^2\,dx+(1-\varepsilon
c_1)\int_{\Omega}|\nabla_xu|^2\,dx\\*+a_0\int_{\Omega_0}|u|^{q+1}\,dx-c_2\varepsilon
\biggl(  \int_{\Omega_0}|u|^{q+1\,dx} \biggr)^{\frac2{q+1}}\leq0.
\end{multline}
From \eqref{4.7} it follows that
\begin{multline}\label{4.10}
\int_{\Omega}|u(x,t)|^{1+q}\,dx\leq(\text{mes}\,\Omega)^{\frac{1-q}2}\biggl(
\int_\Omega |u(x,t)|^2\,dx
\biggr)^{\frac{q+1}2}\\\leq(\text{mes}\,\Omega)^{\frac{1-q}2}\biggl(
\int_\Omega |u_0(x)|^2\,dx \biggr)^{\frac{q+1}2}=
(\text{mes}\,\Omega)^{\frac{1-q}2}y_0^{\frac{q+1}2}=\widetilde{C}=\text{const}
\quad\forall\,t>0.
\end{multline}
Now due to \eqref{4.10} we have
\begin{multline}\label{4.11}
a_0\int_{\Omega_0}|u(x,t)|^{q+1}\,dx-c_2\varepsilon\biggl(
\int_{\Omega_0}|u(x,t)|^{q+1}\,dx
\biggr)^{\frac2{q+1}}\\=\int_\Omega |u(x,t)|^{q+1}\,dx\biggl(
a_0-c_2\varepsilon\biggl(  \int_{\Omega_0}|u(x,t)|^{q+1}\,dx
\biggr)^{\frac{1-q}{1+q}} \biggr)\\\geq\int_\Omega
|u(x,t)|^{q+1}\,dx \bigl(a_0-c_2\varepsilon
\widetilde{C}^{\frac{1-q}{1+q}}\bigr)\geq0
\end{multline}
 if $\varepsilon$ is small enough, namely,
\begin{equation}\label{4.12}
\varepsilon\leq\frac{a_0}{c_2\widetilde{C}^{\frac{1-q}{1+q}}}.
\end{equation}
Thus, if $\varepsilon$ satisfies \eqref{4.12}, then  from \eqref{4.9} it follows
that
$$
\frac d{dt}\int_{\Omega}|u(x,t)|^2\,dx+\varepsilon\int_\Omega
|u(x,t)|^2\,dx\leq0\quad\forall\,t>0.
$$
The last inequality implies the assertion of  Lemma A.3.

{\bf Acknowledgment.} The authors are very grateful to Laurent Veron and Vitali Liskevich
for useful discussions and valuable comments. The second author (AS) has been supported by an INTAS grant through the Program INTAS 05-1000008-7921.

\end{document}